
%
%

\magnification=\magstep1

\font\BBF=msbm10
\def\Bbf#1{\hbox{\BBF #1}}

\font\Small=cmr8
\def\small#1{\hbox{\Small{#1}}}

\font\Big=cmr17

\def\Ker{\hbox {\rm Ker }}

\def\C{\Bbf C}
\def\G{\Bbf G}
\def\P{\Bbf P}
\def\R{\Bbf R}

\def\Z{\Bbf Z}
\def\A{\Bbf A}



\def\qed
{\hskip 10pt
\vrule height 7.5pt depth -0.1pt \vrule height 7.53pt depth -7.2pt width 7.3pt
\hskip -7.5pt \vrule height 0.3pt depth 0pt width 7.6pt \vrule height 7.5pt
depth -0.1pt
}

\def\wt{\widetilde}

\def\sprime{\sp{\prime}}

\def\inv{\sp{-1}}

\def\Im{\hbox{\rm Im }}

\def\Ker{\hbox{\rm Ker }}

\def\Coker{\hbox{\rm Coker }}

\def\pr{\hbox{\rm pr}}
\def\Sing{\hbox{\rm Sing }}

\def\aa{\alpha}
\def\bb{\beta}

\def\gg{\gamma}
\def\ll{\lambda}
\def\mm{\mu}
\def\nn{\nu}
\def\rr{\rho}
\def\ss{\sigma}
\def\ww{\omega}

\def\zz{\zeta}
\def\yy{\eta}

\def\tt{\tau}

\def\vee{\varepsilon}
\def\tthh{\theta}

\def\DD{\Delta}

\def\SS{\Sigma}

\def\medn{\medskip\noindent}
\def\bign{\bigskip\noindent}
\def\parn{\par\noindent}

\def\lra{\longrightarrow}

\def\hs{\hskip 5pt}

\def\hsp#1{\phantom{#1}}
\def\hsa{\hskip 1pt}

\def\llra{\hs\lra\hs}
\def\longmapstoleft{\longleftarrow\hskip -1pt 
\vrule height 4pt width 0.1pt depth -0.3pt}

\def\maprightsp#1{\smash{\mathop{\longrightarrow}\limits\sp{#1}}}

\def\mapstorightsp#1{\smash{\mathop{\longmapsto}\limits\sp{#1}}}

\def\mapstoleftsb#1{\smash{\mathop{\longmapstoleft}\limits\sb{#1}}}
\def\maprightspsb#1#2{\smash{\mathop{\longrightarrow}\limits\sp{#1}\sb{#2}}}

\def\longmapstoleft{\longleftarrow\hskip -1pt 
\vrule height 4pt width 0.1pt depth -0.3pt}

\def\lmt{\longmapsto}

\def\mult{\sp{\times}}

\def\varmaprightsp#1#2{\smash{\mathop{\hbox to #1 {\rightarrowfill}}
\limits\sp{#2}}}
\def\varmaprightsb#1#2{\smash{\mathop{\hbox to #1 {\rightarrowfill}}
\limits\sb{#2}}}

\def\vertcong{\parallel\wr}

\def\hookdownarrow
{{}\sp\cap \hskip -3.703pt \lower 2pt \hbox{$\downarrow$}}
\def\hookuparrow
{\lower 2pt\hbox{${}\sb{\cup}$}\hskip -3.7pt \lower -1pt\hbox{$\uparrow$}}

\def\setbar{\hs ; \hs}

\def\set#1#2{\{\hs{#1}\setbar{#2}\}}

\def\sethd#1#2#3
{\Bigl\{ \hs{#1}\setbar
{\matrix{\hbox{#2} \hfill\cr\hbox{#3} \hfill }} \Bigr\} }
\def\setht#1#2#3#4
{\biggl\{ \hs{#1}\setbar
{\matrix{\hbox{#2} \hfill\cr\hbox{#3} \hfill\cr\hbox{#4} \hfill}} 
\biggr\} }

\def\ibid{[{\it ibid.}]}

\def\diagram#1{
\def\normalbaselines{\baselineskip20pt\lineskip3pt\lineskiplimit3pt}
\matrix{#1}}

\def\varvarmatrix#1#2#3#4
{\def\normalbaselines{\baselineskip#1\lineskip3pt\lineskiplimit3pt}
\vbox{\vskip #2 \hbox{\hfill$\matrix{#3}$\hfill}\vskip #4}}

\def\varmatrix#1#2
{\def\normalbaselines{\baselineskip#1\lineskip3pt\lineskiplimit3pt}
\matrix{#2}}

\def\varpmatrix#1#2
{\def\normalbaselines{\baselineskip#1\lineskip3pt\lineskiplimit3pt}
\pmatrix{#2}}

\def\varpmatrix#1#2
{\def\normalbaselines{\baselineskip#1\lineskip3pt\lineskiplimit3pt}
\pmatrix{#2}}

\def\bigcases#1
{\biggl\{\, \vcenter{\normalbaselines{\mathsurround=0pt}
\ialign{$##\hfil$&\quad##\hfil\crcr#1\crcr}}\biggr.}

\def\st{\subset}
\def\sm{\setminus}

\def\pione{\pi\sb1}

\def\and{\hbox{\rm and}}

\def\and{\hbox{and}}
\def\for{\hbox{for}}
\def\where{\hbox{where}}

\def\bfclaim#1{\medskip\noindent{\bf #1 }\hs }

\def\oover#1#2{\displaystyle{{#1}\over{#2}}}

\def\locus#1{\{ #1 \}}

\def\Gm{\G \sb m}

\def\trgp{\{1\}}

\def\Pt{\P\sp 2}

\def\xxxiii{\xi\sb 0, \xi\sb 1, \xi\sb 2}

\def\Auv{{\Bbf A}\sp 2 \sb{(u, v)}}
\def\As{{\Bbf A} \sp 1 \sb s}
\def\At{{\Bbf A}  \sp 1\sb t}

\def\Bs{B\hskip -1.54pt s}

\def\uv{\sb{uv}}
\def\s{\sb s}
\def\t{\sb t}

\def\sbee{\sb{\vee}}

\def\pv{p\sb v}

\def\afc{s, t, x, y}

\def\smo{\sm \{0 \}}

\def\Degtyarev{[3]}
\def\Abhyankar{[1]}
\def\BruceGiblin{[2]}
\def\Dimca{[4]}
\def\Nemethi{[5]}
\def\Nori{[6]}
\def\Oka{[7]}
\def\Okapre{[8]}
\def\ShimadaF{[9]}
\def\ShimadaW{[10]}
\def\Zariski{[11]}

\def\sps{\phantom{a}}

\def\DimcaPart#1{[4; #1]}
\def\NoriPart#1{[6; #1]}
\def\ShimadaFPart#1{[9;  #1]}
\def\ShimadaWPart#1{[10;  #1]}

\def\Gpq{G\langle p, q \rangle}
\def\Hqk{H(q; k)}
\def\Hqo{H(q; 1)}

\def\aaz{\bb\sp{\mu}\aa\sp{\nu}}

\def\pmat#1#2#3#4{\left[\matrix{#1 &#2\cr #3 & #4\cr}\right]}

\def\vphi{\varphi}

\def\o{\sb 0}

\topskip 1cm 
\centerline
{\Big   On projective plane curves whose complements have }
\bigskip
\centerline
{\Big   finite non-abelian fundamental groups }
\vskip 1cm 
\centerline
{Ichiro Shimada}
\bigskip
\centerline
{Max-Planck-Institut f\"ur Mathematik}
\vskip 1cm
\bign
{\bf \S 0. Introduction}
\medskip
In this paper, we present  new examples of
singular projective plane curves $C \st \Pt$
defined over the complex number field $\C$ 
such that the topological fundamental group 
$\pione (\Pt \sm C)$
is non-abelian and  finite.
\medskip
The first example of a curve with this property is the three cuspidal quartic 
curve
discovered by Zariski in \Zariski.
Let $C\sb 0 \st \Pt$ be a quartic curve defined over $\C$
whose singular locus consists of $3$ cusps.
It is known \BruceGiblin\  that any three cuspidal quartic is projectively
isomorphic to the curve defined by
$$
x\sp2 y\sp 2 + y\sp 2 z\sp 2 + z\sp 2 x\sp 2 - 2xyz(x+y+z) =0.
$$
The topological fundamental group $\pione (\Pt \sm C\sb 0)$
is isomorphic to the binary $3$-dihedral group
$$
\wt D\sb 3  := \langle \hs \aa, \bb \hs |\hs  \aa\sp 2 = \bb\sp 3 = 
(\aa\bb)\sp 2 \hs \rangle
$$
of order $12$ (cf. \Zariski \sps\DimcaPart{Chapter 4, \S 4}).
In \Abhyankar,
Abhyankar studied the complement of the three cuspidal quartic
over an algebraically closed field $k$ of arbitrary characteristics.
He showed that,
if $\hbox{char } k \not = 2, 3$,
then the tame fundamental group
of the complement is isomorphic to $\wt D\sb 3$.
\medskip
It is only quite recent that 
we come to know other examples exist.
In \Degtyarev,
Degtyarev has found an irreducible curve of degree $5$
with three singular points of type $A\sb 4$
such that the complement has a nonabelian fundamental group 
of order  $320$.
After this,
some infinite series of projective plane curves with this property
have been constructed
by Oka  \Okapre\  and, independently, by the author \ShimadaW.
However,
the curves in these infinite series  
should be considered as offspring 
of the three cuspidal quartic,
because they are obtained  as pull-backs of $C\sb 0$
by certain branched coverings of the plane.
On the other hand,
the curves we will construct 
in this paper
can be regarded as siblings of the three cuspidal quartic,
as is explained below.
\medskip
From now on,
we will work over the complex number field $\C$ exclusively.
The fundamental group always means
the topological fundamental group.
\medskip
Let $q$ and $k$ be positive integers,
and suppose that $q$ is odd and $> 1$.
Let $F(s, t, x)$
be the polynomial
$$
F(s, t, x) \quad := \quad
{{(s\sp 2 x + t\sp q )\sp 2 - (s\sp 2 + t\sp 2 ) \sp q }\over{s\sp 2}} \hs .
$$
(Note that $(s\sp 2 x + t\sp q )\sp 2 - (s\sp 2 + t\sp 2 ) \sp q $
is divisible by $s\sp 2$.)
We put weights to the variables as follows;
$$
\deg s = k, \quad \deg t = k, \quad \deg x = (q-2)k.
$$
Then $F$ is weighted homogeneous of total degree $2(q-1)k$.
Recall that a germ of curve singularity
is called to be of type $A\sb{q-1}$
if it is locally defined by $u\sp 2 - v\sp q = 0 $.
The main result of this paper is as follows:
\bfclaim{Theorem.}
{\sl
Let $S(\xxxiii)$, $T(\xxxiii)$ and $X(\xxxiii)$ 
be general homogeneous polynomials of degree $k$, $k$ and  
$(q-2)k$, respectively.
Let $C(q, k) \st \Pt$ be the curve of degree $2(q-1)k$
defined by
the homogeneous equation
$$
F(\hs S(\xxxiii), \hs T(\xxxiii), \hs X(\xxxiii)\hs )\hs  =\hs 0.
$$
Then the singular locus of $C(q, k)$ consists of 
$(2q-3)k\sp 2$ singular points of type $A\sb{q-1}$,
and the fundamental group $\pione (\Pt \sm C(q, k))$
is a central extension of the dihedral group $D\sb q$ of order $2q$
by the cyclic group $\Z / (k(q-1))$.
In particular, $\pione (\Pt \sm C(q, k))$ is finite and non-abelian.
}
\medn
Actually, we are going to describe 
the structure of $\pione (\Pt \sm C(q, k))$
completely as a subgroup of $GL\sb 2 (\C)$.
\medskip
The  curve $C(3, 1)$ is the three cuspidal quartic.
The curves $C(q, k)$ $(k \ge 2)$
can be  obtained from $C(q, 1)$
by the method of \ShimadaW.
So the essentially new part of Theorem is the construction of $C(q, 1)$
with $q$ odd integer $\ge 5$. 
\medskip
The idea of  the construction
stems from the following observation due to Zariski \Zariski.
The three  cuspidal quartic curve
is an irreducible component
of a reducible  sextic curve defined by  an equation of the form
$\phi\sb 2 \sp 3 + \phi \sb 3\sp 2 = 0 $, where
$\phi\sb{\nn}$ is a suitable homogeneous polynomial on $\Pt$  of degree $\nn$,
with the residual component being a double line.
\medskip
Let $p$, $q$ be integers $>1$ such that $(p, q)=1$,
and $k$ a positive integer.
Let $f$ and $g$ be general homogeneous polynomials on $\Pt$  of degree $qk$ 
and $pk$,
respectively, 
and let $C\sb{f, g}$ be the curve defined by $f\sp p - g\sp q = 0 $.
We denote by  $G$  the group 
$\langle \aa, \bb  |   \aa\sp p = \bb\sp q \rangle$.
In \ShimadaW, it was shown that
the fundamental group $\pione (\Pt \sm C\sb{f, g})$  is isomorphic to
the group
$$
\langle \hs a, b, c \hs | \hs a\sp p=b\sp q=c, \hs c\sp k = e \hs \rangle,
$$
which is obtained by putting one more relation $\aa\sp{pk} = e$ on $G$.
In particular, when $k=1$, this group is isomorphic to the free product 
$\Z/(p) * \Z/(q)$
of cyclic groups of order $p$ and $q$ (\Oka, \Nemethi).
The singular locus of the curve $C\sb{f, g}$ consists of
$pqk\sp 2 $ points at which $C\sb{f, g}$ 
is locally defined by $u\sp p + v\sp q = 0$.
Inspired by the above observation of Zariski,
we will consider a situation
in which the curve $C\sb{f, g}$ degenerates 
into the union
of a reduced curve $C\sb 1$ and a non-reduced curve  $mC\sb 2 $
with  multiplicity $m \ge 2$.
More precisely, we specialize $f$ and $g$
to $S\sp m X + T\sp q$ and $S\sp m Y + T\sp p$,
respectively,
where $S$, $T$, $X$ and $Y$ are general homogeneous polynomials on 
$\Pt$ of suitable 
degrees.
It turns out that,
in many cases,
$\pione (\Pt \sm C\sb 1)$ remains non-abelian.
Moreover,
there is a surjection from $G$ to $\pione (\Pt \sm C\sb 1)$,
but the generators $\aa$    and $\bb$ acquire  more relations than
$\aa\sp{pk}= e$ in $\pione (\Pt\sm C\sb 1)$.
In the particular case described in Theorem --
that is, when $p=2$ and $m=2$ --
the added  relations are so imposing that the quotient 
$\pione (\Pt \sm C\sb 1)$
becomes finite.
Moreover, the curve $C\sb 1$ does not acquire any other type of 
singularities than those of $C\sb{f, g}$.
\medskip
The main tools of the calculation of $\pione (\Pt \sm C(q,k))$
are the weighted Zariski's hyperplane section theorem \ShimadaWPart{Theorem 1},
and the generalized Zariski-Van Kampen's theorem \ShimadaFPart{Theorem 2}.
The main theorem stated above is 
obtained by combining Propositions 2 and 3 in sections 3 and 4.
\medskip
{\bf Acknowledgment.}
\hs
The problem to search for the singular projective plane curves
whose complements have non-abelian and finite fundamental groups
was suggested by Professor M.\ Oka
at the workshop on ``Fundamental groups and branched coverings''
held at Tokyo Institute of Technology on December 1994.
The author would like to thank him for stimulating discussions.
\par
Part of this work was done during the author's stay at 
Max-Planck-Institut f\"ur Mathematik.
He thanks to people at this institute for their warm hospitality.
\bign
{\sl Conventions}
\medskip
(i)
Let $\aa$ and $\bb$ be  paths in a topological space.
We define the product of $\aa$ and $\bb$
in such a way
that $\aa\bb$
is defined if and only if the ending point of $\aa$ coincides with 
the starting point of $\bb$.
\par
(ii)
In order to distinguish various affine spaces,
which will appear during the course of the discussions,
we write the affine coordinate(s)
of the space at the bottom of the name of the space.
For example, $\A\sp 1\sb u$ denotes 
an affine line with an affine coordinate $u$ (that is,  a $u$-line),
and $\Auv$ is the affine plane $\A\sp 1 \sb u \times \A\sp 1 \sb v$.
\par
(iii)
A loop in a topological  space is represented by a continuous map from
the circle 
$\{ e\sp{i\tthh} \in \C | \tthh \in \R \}$
with the base point $1$ and the counter-clockwise orientation.
For example,
instead of  denoting a loop $[0,1] \to \Auv$ in such a way as 
$t\mapsto (u, v) = ( e\sp{2\pi i m t}, e\sp{2\pi i n t} )$,
we simply write {\sl a loop $(e\sp{im\tthh}, e\sp{in\tthh})$ on $\Auv$.}
\bign
{\bf \S 1. The defining polynomial of the curve}
\medskip
We consider a more general situation 
than the one described in Theorem.
\medskip
Let $F(s, t, x, y)$
be the polynomial
$$
F(s, t, x, y) 
\quad = \quad 
{{(s\sp m x + t\sp q )\sp p - (s\sp m y + t\sp p)\sp q}\over{s\sp m}},
\eqno{(1.1)}
$$
where $p, q$ and $m$ are integers $> 1$
such that $(p, q) =1$.
Let $k$ and $l$ be positive  integers such that $pl \ge  mk$ and $ql \ge mk$.
We put weights to the variables as follows;
$$
\deg s = k, \quad \deg t = l, 
\quad \deg x = ql -mk \quad\and \quad \deg y = pl-mk.
\eqno{(1.2)}
$$
Then the polynomial $F$ becomes 
weighted homogeneous of total degree $pql - mk$.
Let $S(\xxxiii)$, $T(\xxxiii)$, $X(\xxxiii)$ and $Y(\xxxiii)$
be {\sl general} homogeneous polynomials of 
degree $k$, $l$, $ql-mk$ and $pl-mk$,
respectively.
Substituting $s,t,x,y$ with $S, T, X, Y$,
we get a homogeneous polynomial $F(S, T, X, Y)$
of degree $pql-mk$ in variables $(\xxxiii)$.
Therefore, we get a projective plane curve 
$$
C \quad := \quad \{  F(S, T, X, Y) = 0 \} \quad \st \quad \Pt
$$
of degree $pql - mk$.
This $C$ is the main object of the investigation in this article.
We shall give a group presentation of
the fundamental group $\pione (\Pt \sm C)$ in section 2.
In section 3, we shall investigate the structure of this group under
the assumptions  $p=m=2$ and $k=l$.
In section 4, we shall study the singularities  of $C$ 
under the assumption  $p=m=2$.
\bign
{\bf \S 2. Group presentation of $\pione (\Pt\sm C)$}
\medskip
Let $U \st \A\sp 4\sb{(\afc)}$ be the complement of the hypersurface
$\SS := \{ F(\afc) = 0 \}$.
Then the multiplicative group $\G\sb m $ 
of non-zero complex numbers acts on $U$
by the weights (1.2); that is, 
$$
\ll \cdot ( \afc ) =
 (\hs \ll \sp{k} s\hs ,\hs  \ll\sp{l} t\hs ,\hs  \ll\sp{ql-mk} x\hs ,
\hs  \ll\sp{pl-mk} y\hs)
\quad\where\quad \ll \in \G\sb m, \eqno{(2.1)}
$$
because $F$ is weighted homogeneous under the weights (1.2).
Let 
$$
\rr \quad : \quad \pione (\G\sb m) \llra \pione (U)
$$
be the natural homomorphism induced by this action.
Let $\mm$ and $\nn$ be 
integers such that $\mm p + \nn q = 1$.
(Recall that $p$ and $q$ are supposed to be prime to each other.)
\bfclaim{Proposition 1.}
{\sl
The fundamental group  $\pione (U)$ is isomorphic to
$$
\langle \hs \aa, \hs  \bb \hs | \hs \aa\sp p = \bb\sp q, \hs 
\aa (\aaz)\sp m = (\aaz)\sp m \aa,  
\hs \bb(\aaz)\sp m = (\aaz)\sp m \bb \hs \rangle,
$$
and, under this isomorphism,  the image of 
the counter-clockwise generator of $\pione (\G\sb m)$
by $\rr$
is given by $\bb\sp{ql} (\aaz)\sp{-mk}$ 
in terms of the generators $\aa$ and $\bb$.
}
\medn
{\it Proof.}
\hs
The plan of the proof is as follows.
The first projection from $U$ to the $s$-line $\As$
is a locally trivial fiber space over $\As \smo$.
Moreover,
a general fiber is homotopically equivalent
to the space
$$
\Auv \sm E, \quad\hbox{where $E$ is the affine plane 
curve defined by $u\sp p - v \sp q = 0$},
$$
and hence its fundamental group is isomorphic to
$\langle \aa, \bb | \aa\sp p = \bb\sp q \rangle$ 
(cf. \DimcaPart{Chapter 4, \S 2}).
The kernel of the surjection
$\langle \aa, \bb | \aa\sp p = \bb\sp q \rangle \to \pione (U)$
induced from the inclusion of a general  fiber
into the total space $U$
can be determined by \ShimadaFPart{Theorem 2}.
\medskip
Let $\pr\sb 1 : U\to \As$ be the first projection.
The fiber $\pr\sb1\inv (\ss)$ over the point $s = \ss$
will be denoted by $U\sb{\ss}$.
We put 
$$
U\sp o \quad := \quad U\sm U\sb 0 .
$$
Over $\As \smo$,
the affine hypersurface $\SS$ is defined
by $(s\sp m x + t\sp q)\sp p - (s\sp m y + t\sp p)\sp q = 0$.
Hence there exists an isomorphism $\Phi$
from $U\sp o$ to $(\As\smo)\times\At\times(\Auv\sm E)$
which makes the following diagram commutative;
$$
\diagram{
U\sp o \hskip -20pt&& \quad\maprightspsb{\sim}{\Phi}\quad && (\As\smo ) 
\times \At \times (\Auv \sm E) \cr
& \pr\sb 1 \hs\searrow && \swarrow &\hskip -50pt
\hbox{\small{the first projection}}\cr
&& \As\smo. &&\cr
}
\eqno{(2.2)}
$$
This $\Phi$ is given by
$$
\varmatrix{16pt}{
(s, t, x, y) & \quad\mapstorightsp{\Phi}\quad &
(s, t, \hs \hs(s\sp m x + t\sp q, s\sp m y + t\sp p)) \cr
(s, t, \hs {{(u-t\sp q)}/{s\sp m}}, \hs {{(v - t\sp p) }/{s\sp m}})
&\quad\mapstoleftsb{\Phi\inv} \quad & \hskip -30pt (s, t, \hs (u, v)). \cr
}
$$
We put
$$
\Phi =( \phi\s, \phi\t, \phi\uv ), 
\qquad\where\quad \phi\s : U\sp o \to \As\smo, 
\quad\phi\t : U\sp o \to \At, \quad  \phi\uv : U\sp o \to \Auv \sm E.
$$
Then we have an isomorphism
$$
\pione (U\sp o , b) \quad \cong \quad 
\pione (\As\smo,\hs  \phi\s (b)) \hs\times\hs 
\pione (\Auv \sm E, \hs \phi\uv (b)) \eqno{(2.3)}
$$
induced by $\Phi$.
As a base point $b$,
we choose 
$$
b = (\vee, \hs 1,\hs  1/\vee,\hs  0) \quad \in U\sb{\vee}\st U\sp o \st U
$$
where $\vee$ is a small positive real number.
Then we have
$$
\phi\s (b) = \vee \quad \and  \quad\phi\uv (b) =(1+\vee\sp{m-1}, 1).
$$
We shall write down explicitly the isomorphisms
$$
\eqalign{
\pione (\As\smo, \hs \phi\sb s (b)) \hsp{aa}&\quad\cong\quad \Z, 
\qquad \and \cr
\pione (\Auv \sm E, \hs \phi\uv (b) ) &\quad\cong\quad \Gpq, 
}
$$ 
where
$$ 
\Gpq := \langle\hs \aa, \bb \hs |\hs \aa\sp p = \bb\sp q \hs \rangle.
$$
\par
As usual,
for $\pione (\As\smo, \hs \phi\sb{s} (b))\cong\Z$,
we take the homotopy equivalence class of the loop
$s=\vee e\sp{i\theta}$ on $\As\smo$ as
a positive generator $1 \in \Z$.
\medskip
For the second isomorphism,
whose proof can be found
in \DimcaPart{Chapter 4, \S 2},
we consider the second projection $\pv :\Auv \sm E \to \A\sp 1\sb v$.
The fiber $\pv\inv (1)$ is,
by the first projection,
identified with the $u$-line $\A\sp 1 \sb u$ minus the $p$-th roots of unity.
For each integer $j$, we put
$$
a\sb j\sprime  \quad 
:= \quad [\ww(j)\yy (j)\ww(j)\inv] \quad \in \quad
\pione (\pv\inv (1), \phi\uv (b)),
$$ 
where $\ww(j) : [0, 1] \to \pv\inv (1) $ 
and $\yy(j) : [0, 1] \to \pv\inv (1)$ 
are the paths given below. 
(Here $\pv\inv (1)$ is regarded as a Zariski open subset of the $u$-line.)
$$
\varmatrix{20pt}{
\ww(j) &\quad :\quad&t\sb 1 &\quad\lmt\quad&
u= (1+\vee\sp{m-1} )\exp  (2\pi i  j t\sb 1 / p), \hfill\cr
\yy(j) &\quad :\quad&t\sb 2&\quad\lmt\quad&
u= (1 + \exp (2\pi i  t\sb 2)  \vee\sp{m-1})\exp(2\pi i  j  / p).\cr
}
$$
These $a\sb j\sprime$ generate $\pione (\pv\inv (1), \phi\uv (b))$.
Then the inclusion
$\pv\inv (1) \hookrightarrow \Auv \sm E$
induces a surjective homomorphism on the fundamental groups by 
\NoriPart{Lemma 1.5(C)},
because $\A\sp 1 \sb v$ is simply connected.
We denote  by $a\sb j \in \pione (\Auv \sm E, \phi\uv (b))$
the image of $a\sb j\sprime  \in \pione (\pv\inv (1), \phi\uv (b)) $.
We put
$$
\eqalign
{
\aa&:=  a\sb {q-1} a\sb {q-2} \cdots \cdots a\sb 1 a\sb 0, \quad \and  \cr
\bb&:=  a\sb {p-1} a\sb {p-2} \cdots  a\sb 1 a\sb 0.
 \cr
}
$$ 
Then we get the identification 
$$
\eqalign{
\pione (\Auv \sm E, \phi\uv (b)) 
& = \langle a\sb j \hs (j\in \Z)
\hs  | \hs a\sb{j+q} = a\sb j, \hs a\sb{j+p} = \bb a\sb j \bb\inv 
\hs\hbox{  for all $j \in \Z$ }\hs \rangle\cr
& = \langle \aa, \bb \hs  | \hs  \aa \sp p = \bb \sp q \rangle = \Gpq, 
}
$$
(cf. \DimcaPart{Chapter 4, \S 2}).
Note that we have 
$$
a\sb 0  =  \aaz, \qquad \where\quad \mm p + \nn q = 1, \eqno{(2.4)}
$$
(cf. \ibid).
Combining these isomorphisms and (2.3),
we have fixed an isomorphism
$$
\pione (U\sp o , b ) \quad \cong \quad \Z \hs \times \hs \Gpq .
\eqno{(2.5)}
$$
\par
Consider the fiber $\pr\sb 1 \inv (\vee) = U\sbee$,
which contains the base point $b$.
By (2.2) and  (2.3),
the image of the natural homomorphism $\pione (U\sbee, b) \to 
\pione (U\sp o , b)$
induced by the inclusion 
can be identified 
with the subgroup $\{ 0 \} \times \Gpq$ 
of $\Z \times \Gpq$
via (2.5).
Note that, by \NoriPart{Lemma 1.5(C)},
the natural homomorphism $\pione (U\sbee, b) \cong \Gpq\to \pione (U, b)$
induced by the inclusion 
is surjective, because $\As$ is simply connected.
We shall determine the kernel of this surjection.
\bfclaim{Claim 1.}
{\sl The open disk
$$
D = \set{ (z, 1, 1/\vee, 0) }{ |z|< 2 \vee } \quad\st\quad \A\sp 4\sb{(\afc)}
$$
is contained in $U$.
}
\medn
{\it Proof.}
\hs
On the hyperplane defined by 
$s=0$, the hypersurface $\SS \cap \{ s = 0 \} \st \A\sp 3 \sb{(t, x, y)}$
is defined by
$pxt\sp{q(p-1)} - q y t \sp{p(q-1)} = 0.$
The intersection point $ (0, 1, 1/\vee, 0)$ of $D$ and $\{ s=0\}$
does not satisfy this equation.
When
$s\not = 0$,
the hypersurface $\SS$ is defined by
$(s\sp m x + t \sp p ) \sp q - ( s \sp m y + t \sp q ) \sp p = 0$. 
Putting  $s=z$, $t=1$, $x = 1/ \vee$ and $y=0$, we get
$(z\sp m / \vee + 1 ) \sp q - 1 = 0$.
Since  $m \ge 2$,
this equation does not have any roots in $0 < |z| < 2 \vee$,
because $\vee$ is small enough.
Hence we get $D\cap \SS =\emptyset$. \qed
\bfclaim{Claim 2.}
{\sl
The kernel of the surjective homomorphism
$\pione (U\sbee, b) \cong \Gpq \to \pione (U, b)$ is generated by the set 
$$
\set{\gg \inv a\sb 0 \sp{-m} \gg a\sb 0 \sp{m} }{\gg \in \Gpq },
$$
where $a\sb 0 = \bb\sp{\mm}\aa\sp{\nn}$ by (2.4).
}
\medn
In other words,
$\pione (U, b)$ is isomorphic to the maximal quotient
group of $\Gpq$
such that the image of $a\sb 0 \sp m$
is contained in the center.
\medn
By this claim, the first assertion of Proposition 1 will be proved.
\medn
{\it Proof.}
\hs
Let $\DD := \set{s}{|s| < 2\vee} \st \As$
be the open disk on the $s$-line with the origin $0$ and the radius $2\vee$.
Let $U\sb{\DD}$ be $\pr\sb 1 \inv (\DD)$.
Since $\pr\sb 1 : U \to \As$ is a locally trivial fiber space 
over $\As\smo$, $U\sb{\DD}$ is a strong deformation retract of $U$.
Therefore
the inclusion $U\sb{\DD} \hookrightarrow U$
induces an isomorphism
$\pione (U\sb{\DD}, b) \cong \pione (U, b)$.
Thus it is enough to consider $\Ker (\pione (U\sbee, b) \to 
\pione (U\sb{\DD}, b))$.
\medskip
Note that, by Claim 1, the projection $\pr\sb 1 : U \to \As$ has a section
$$
\ss \hs : \hs s \hs \lmt \hs (s, 1, 1/\vee, 0)
$$
over $\DD$ passing through the base point $b$.
We put $\DD\mult := \DD\smo$.
By the section $\ss$,
the group $\pione (\DD\mult, \vee)$ acts on $\pione (U\sbee, b)$
in a natural way,
because $\pr\sb 1 : \pr\sb 1 \inv (\DD\mult) \to \DD\mult$ 
is a locally trivial fiber space.
Let
$\zz \in \pione (\DD\mult, \vee)$ be the counter-clockwise generator
of $\pione (\DD\mult, \vee)\cong \Z$,
and let $\gg \mapsto \zz\sb* ( \gg)$ denote 
its monodromy action on $\gg \in \pione (U\sbee, b)$.
By \ShimadaFPart{Theorem 2},
the kernel of $\pione (U\sbee, b) \to \pione (U\sb{\DD}, b)$
is generated by the set
$$
\set{\gg\inv \cdot \zz\sb* (\gg) }{ \gg \in \pione (U\sbee, b)}.
$$
There is an isomorphism $\Psi$ from $\pr\sb 1 \inv (\DD\mult)$
to $\DD\mult \times U\sb{\vee}$
which makes the following diagram commutative;
$$
\diagram{
\pr\sb 1 \inv (\DD\mult)  \hskip -20pt&& 
\quad\maprightspsb{\sim}{\Psi}\quad && 
\DD\mult \times U\sbee \cr
& \pr\sb 1 \hs\searrow && \swarrow &\hbox{\small{the first projection}}\cr
&& \DD\mult. &&\cr
}
$$
This $\Psi$ is given by
$$
\varmatrix{16pt}{
(s, \hs t, \hs  x, \hs y) & \quad\mapstorightsp{\Psi}\quad &
(s, \hs (\vee,\hs  t, \hs s\sp m x / \vee\sp m, \hs s\sp m y / \vee\sp m)) \cr
(s, \hs t, \hs \vee\sp m x\sprime / s\sp m, \hs \vee\sp m y\sprime / s\sp m)
&\quad\mapstoleftsb{\Psi\inv} \quad & 
\hskip -30pt (s, \hs (\vee, \hs t, \hs x\sprime, \hs y\sprime)). \cr
}
$$
Hence the fiber structure  of  $ \pr\sb 1 : 
\pr\sb 1 \inv (\DD\mult ) \to \DD\mult$
is trivial. 
However
the motion of the base point causes a non-trivial monodromy action
of $\pione (\DD\mult)$ on $\pione (U\sb{\vee})$.
We let the generator  $\zz \in \pione (\DD\mult, \vee)$ be represented by
the loop $\vee e\sp{i\tthh}$ on $\DD\mult$.
The image of this loop by the morphism 
$$
\DD\mult \quad \maprightsp{\ss} \quad \pr\sb 1 \inv (\DD\mult ) 
\quad  \maprightsp{\Psi} \quad
\DD\mult \times U\sbee \quad \maprightsp{\pr\sb 2 } \quad U\sbee
$$
is the loop $(\vee, 1, e\sp{im\tthh}/\vee, 0)$
from $b$ to $b$.
Let $\wt \zz \in \pione (U\sbee, b)$ be the homotopy equivalence class of 
this loop.
By Lemma below,
the action $\gg\mapsto \zz\sb*(\gg)$ is given by
$\zz\sb * (\gg) = \wt\zz \inv \cdot \gg \cdot \wt \zz$.
We shall write down 
the image of $\wt \zz$ by the isomorphism $\pione (U\sbee, b) \cong \Gpq$.
The image of the loop $(\vee, 1, e\sp{im\tthh}/\vee, 0)$ on $U\sbee$
by $\phi\uv | \sb{U\sbee} : U\sbee \to \Auv \sm E$
is the loop  $(1 + \vee\sp{m-1} e\sp{im\tthh}, 1)$ on $\Auv \sm E$,
which is nothing but the loop $(\yy(0))\sp m$.
This  represents $a\sb 0\sp m \in \Gpq$.
Thus we get
$$
\zz\sb * (\gg) =  a\sb 0\sp{-m} \gg a\sb 0\sp m,
$$ 
and the proof of Claim 2 is completed. \qed
\bfclaim{Lemma.}
{\sl Let $(V, b\sb V)$ and $(W, b\sb W)$
be topological spaces with base points.
Suppose that we are given
a section
$\ss : V \to V\times W$
of the first projection $V\times W \to V$
such that $\ss (b\sb V) = (b\sb V, b\sb W)$.
Then the action of $\aa \in \pione (V, b\sb V)$ on $\gg \in \pione (W, b\sb W)$
induced by this section is
given by 
$\gg \mapsto \wt\aa\inv\cdot \gg \cdot \wt \aa$,
where $\wt \aa$ is the image of $\aa $
by the homomorphism  $\pione (V, b\sb V) \to \pione (W, b\sb W)$
induced from the composition of the section
$\ss : V\to V\times W$ and the second projection
$\pr\sb 2 : V\times W \to W$.
}
\qed
\medskip
Next, we will investigate the image of $\rr : \pione (\G\sb m) \to \pione (U)$
induced by the action (2.1).
By Claim 1,
the loop 
$(\vee e\sp{i\tthh}, 1, 1/\vee, 0)$
on $U\sp o$ represents an element of the kernel
of the homomorphism  $\pione (U\sp o, b) \to \pione (U, b)$
induced by the inclusion.
By $\phi\uv : U\sp o \to \Auv \sm E$,
this loop is mapped to the loop
$(1 + \vee\sp{m-1} e\sp{im\tthh}, 1)$,
which represents $a\sb 0 \sp m \in \pione (\Auv \sm E, \phi\uv(b)) \cong \Gpq$.
The image of the class of this loop by 
$\phi\sb{s*}: \pione (U\sp o, b) \to \pione(\As\smo, \phi\sb s (b))$
is obviously the generator $1 \in \pione(\As\smo, \vee) \cong \Z$.
Thus the loop 
$(\vee e\sp{i\tthh}, 1, 1/\vee, 0)$
represents $(1, a\sb 0\sp m) \in \Z \times \Gpq$
via the isomorphism (2.5).
Hence we have shown that
$$
(1, a\sb 0\sp m ) \in  \Ker (\pione (U\sp o, b) \to \pione (U, b)).
\eqno{(2.6)}
$$
Since $U\sp o \st U$ is invariant under the action (2.1) of $\G\sb m$,
there is a natural homomorphism $\pione (\G\sb m ) \to \pione (U\sp o)$,
which factors $\rr$.
In order to determine the image of
$\pione (\G\sb m ) \to \pione (U\sp o)\cong \Z \times \Gpq$,
we choose another base point
$$
c := (1,0,R, 0 ) \in U\sp o,
$$
where $R$ is a positive real number large enough.
\medskip
{\it Remark.}\hs
A path from $c$ to $b$ in $U\sp o$
induces an isomorphism $\pione (U\sp o, c) \cong \pione (U\sp o, b)$,
which depends on the homotopy class of the path.
However,
the image of $1 \in \Z \cong \pione (\G\sb m)$
in $\pione (U\sp o, c)$
always corresponds to the image of $1$
in $\pione (U\sp o, b)$,
whatever choice we may have made on the path,
because the image of $\pione (\G\sb m) \to \pione (U\sp o)$ is 
contained in the center.
\medskip
The orbit of the new base point $c=(1,0,R,0)$ by the loop $e\sp{i\tthh}$ 
in $\G\sb m$
is  the loop  
$(e\sp{ik\tthh}, 0, R e\sp{i(ql-mk)\tthh}, 0)$ in $U\sp o$,
which is mapped by $\phi\s$ and $\phi\uv$
to the loops
$e\sp{ik\tthh}$ in $\As\smo$ and 
$(R e\sp{iql\tthh}, 0)$ 
in $\Auv\sm E$,
respectively.
The loop $e\sp{ik\tthh}$ on $\As\smo$ represents
$k \in \Z \cong \pione (\As\smo)$.
Note that, by the definition, the element $\bb \in \Gpq$ is,
as an element of $\pione (\Auv \sm E)$, 
represented by
the loop 
$(e\sp{i\theta} (1+ \vee\sp{m-1}), 1)$ 
on $p\sb v \inv (1)$.
The loop $(R e\sp{iql\tthh}, 0)$ 
on $\Auv\sm E$, which is supported on $\pv\inv (0)$,
can be deformed to the loop
$(e\sp{iql\theta} (1+ \vee\sp{m-1}), 1)$ 
on $\pv\inv(1)$,
which represents $\bb\sp{ql} $ in $\pione (\Auv\sm E) \cong \Gpq$.
Therefore, by Remark above,
the image of $1 \in \Z \cong \pione (\G\sb m)$
in $\pione (U\sp o, b)$ 
is 
$$
(k, \bb\sp{ql}) \in \Z \times \Gpq\cong \pione (U\sp o, b).
$$
Because of (2.6),
the elements $(k, \bb\sp{ql})$ and $(0, \bb\sp{ql}a\sb 0\sp{-mk})$
are mapped to the same element in $\pione (U)$.
Therefore, when we regard $\pione (U)$
as a quotient of $\pione(U\sbee) \cong \Gpq$,
the image of $1 \in \Z \cong \pione (\G\sb m)$ in $\pione (U)$ by $\rr$
is the image of 
$\bb\sp{ql}a\sb 0\sp{-mk} \in \Gpq$.
\qed
\medn
\bfclaim{Corollary.}
{\sl
The fundamental group $\pione (\Pt\sm C)$
is isomorphic to the group
$$
\biggl\langle \aa,\bb \hs  \biggr|\hs  
\varmatrix{18pt}{
\aa\sp p =\bb\sp q, & \hs \bb\sp{ql} = (\bb\sp{\mm}\aa\sp{\nn})\sp{mk} \cr
\aa(\bb\sp{\mm}\aa\sp{\nn})\sp m = (\bb\sp{\mm}\aa\sp{\nn})\sp m \aa, &
\bb(\bb\sp{\mm}\aa\sp{\nn})\sp m = (\bb\sp{\mm}\aa\sp{\nn})\sp m \bb } 
\hs  \biggr\rangle.
$$
} 
\medn
{\it Proof.}
\hs
Note that  the affine hypersurface $\SS$ is reduced.
Therefore, 
by \ShimadaWPart{Theorem 1},
when both of $\deg x$ and $\deg y$ are positive,
$\pione (\Pt \sm C)$ is isomorphic to the cokernel 
of the natural homomorphism   $\rr : \pione (\G\sb m) \to \pione (U)$,
and thus we get the  group presentation
of $\pione (\Pt\sm C)$ above from Proposition 1.
\medskip
In the case when $\deg x= pl-mk=0$ or $\deg y = ql-mk=0$,
we cannot conclude $\pione (\Pt\sm C) \cong \Coker \rr$
by applying \ShimadaWPart{Theorem 1} directly.
However, this isomorphism still holds in this case, as is shown below.
Interchanging $p$ and $q$ if necessary,
we may assume that $pl-mk >0$ and $ql-mk=0$.
Let $\pr\sb 4 : U\to \A\sp 1\sb y $ be the projection
to the $y$-line.
If $\yy \not = 0$,
then the fiber $\pr\sb 4 \inv (\yy)$ is isomorphic to
$\pr\sb 4 \inv (1)$.
Indeed,
the morphism
$(s, t, x, 1) \mapsto (s\yy\sp{-1/m}, t, x \yy, \yy)$
gives an isomorphism from $\pr\sb 4 \inv (1)$ 
to $\pr\sb 4 \inv (\yy)$.
Hence $\pr\sb 4: U \to \A\sp 1 \sb y$ is a locally trivial fiber space 
over $\A\sp 1 \sb y \smo$.
On the other hand,
it is easy to see that the hypersurface
$\SS \cap \locus{y=0} = \locus{F(s,t,x,0)=0}$ in 
$\A\sp 3 \sb{(s, t, x)}$ is reduced.
Therefore,
by \ShimadaFPart{Theorem 1},
we see that the inclusion $\pr\sb 4 \inv (\yy) \hookrightarrow U$
induces an isomorphism 
$\pione (\pr\sb 4 \inv (\yy) ) \cong \pione (U)$ for $\yy \not = 0$.
On the other hand,
the action (2.1) of $\Gm$ on $U$ leaves $\pr\sb 4 \inv (\yy)$ invariant
because of $\deg y = 0$.
Hence we have a commutative diagram
$$
\varmatrix{14pt}{
& &&\pione (\pr\sb 4 \inv (\yy))\cr
&\nearrow&&\cr
\pione (\Gm) \hsp{aa} && &\vertcong \cr
&\searrow&& \cr
& & & \pione (U). }
\qquad \for \qquad\yy\not= 0
\eqno{(2.7)}
$$
Suppose that the polynomial $Y$, which is  of degree $0$ in this case,
is a constant $\yy\sb 0 \not = 0$.
By \ShimadaWPart{Theorem 1},
$\pione (\Pt \sm C)$ is isomorphic to 
$\Coker (\pione (\Gm) \to \pione (\pr\sb 4 \inv (\yy\sb 0)))$
because of $\deg s > 0$, $\deg t > 0$ and $\deg x >0$.
Hence, by  the diagram (2.7),  we obtain 
$\pione (\Pt \sm C) \cong \Coker(\pione (\Gm) \to \pione (U))$.
\qed
\bign
{\bf \S 3. Structure of $\pione (\Pt \sm C)$ in the dihedral case}
\medskip
The dihedral case in the title of this section
is the case when $p=m= 2$.
In this case,
$q$ is an odd integer $q=2r + 1 $ with $r \ge 1$.
\bfclaim{Proposition 2.}
{\sl Suppose that $p=m=2$,
and $k=l$.
Then the fundamental group $\pione (\Pt \sm C)$
is isomorphic to
a central extension of the dihedral group $D\sb q$
of order $2q$ by the cyclic group $\Z / (k(q-1))$.
In particular, $\pione (\Pt \sm C)$ is non-abelian and 
finite.
}
\medskip
In this section,
we use the following notation;
$$
e(\aa) := \exp(2\pi i \aa).
$$
We consider the dihedral group $D\sb q$ 
as a subgroup of $PGL\sb 2 (\C)$ generated by
$$
\pmat{0}{1}{1}{0}
\quad\and\quad
\pmat{e(1/q)}{0}{0}{1}.
$$ 
{\it Proof.}
\hs
Note that, in the case $p=2$ and $q=2r+1$,
we have $a\sb 0 = \bb\sp{-r}\aa$ by (2.4).
Therefore, by Corollary of Proposition 1,
the group $\pione (\Pt\sm C)$ is isomorphic to
$$
H(q;k) :=
\biggl\langle \aa,\bb \hs  \biggr|\hs  
\varmatrix{18pt}{
\aa\sp 2 =\bb\sp q, & \hs \bb\sp{qk} = (\bb\sp{-r} \aa)\sp{2k} 
\hs \hs & \cdots \hbox{Relations (1)}\cr
\aa(\bb\sp{-r} \aa)\sp 2 = (\bb\sp{-r} \aa)\sp 2 \aa, &
\bb(\bb\sp{-r} \aa)\sp 2 = (\bb\sp{-r} \aa)\sp 2 \bb & \cdots 
\hbox{Relations (2)}} 
\hs  \biggr\rangle.
$$ 
Note that when $k=1$, the relations (2) follows automatically from the 
relations (1).
\bfclaim{Claim.} 
{\sl The order of $ H(q;k) $ is at most $\ 2q(q-1)k$.}
\medn
{\it Proof.}
\hs
Note that in the present situation, we have $\deg y =0$.
Let $\rr\sb{\nn} : \pione (\Gm) \to \pione (U)$
be the homomorphism induced by the action of $\Gm$
on $U$
with weights $(\nn, \nn, (q-2)\nn, 0)$ on the variables $(s, t, x, y)$.
Then  $\Hqk \cong \Coker \rr\sb k$ holds for all $k\ge 1$,
as is shown in the proof of Corollary.
We see that  $\rr\sb k = \rr\sb 1 \circ \ss\sb k$,
where $\ss\sb k : \pione (\Gm) \to \pione (\Gm)$
is induced from the morphism $\ll\mapsto\ll\sp k$.
Hence we have an exact sequence
$$
\Coker \ss\sb k \quad \lra \quad 
\Coker \rr\sb k \quad \lra \quad
\Coker \rr\sb 1 \quad \lra \quad \trgp.
$$
Consequently the order of $\Hqk$ is at most $k$ times the order of $\Hqo$.
Therefore, it is enough to show that the order of $\Hqo$ is at most $2q(q-1)$.
In $\Hqo$,
we have  $\aa\sp 2 = \bb\sp q = (\bb\sp{-r}\aa)\sp 2$,
and hence 
$$
\bb\sp r \aa = \aa \bb \sp{-r}. \eqno{(3.1)}
$$
Any integer $n$ can be written as $nq-2nr$.
Therefore we have 
$$
\bb\sp n \aa = (\bb \sp q)\sp n (\bb \sp r) \sp{-2n} \aa
= (\aa\sp 2 )\sp n \aa (\bb\sp{-r} )\sp{-2n}
= \aa\sp{2n +1} \bb\sp{2nr}.
$$
We also have $\bb\sp n \aa\sp{-1} = \aa\sp{2n-1}\bb\sp{2nr}$.
By repeating these transformations,
every word in the letters $\aa\sp{\pm 1}$ and $\bb\sp{\pm 1}$
can be transformed into a word of the type $\aa\sp M \bb \sp N$.
Since $\aa\sp 2 = \bb\sp q$, we can assume that $M=0$ or $1$.
Therefore the order of $\Hqo$ is at most $2$ times the order of 
$\bb$ in $\Hqo$.
From the relation $\bb\sp{-r} \aa\bb\sp{-r}\aa = \aa\sp 2 = \bb\sp{2r +1}$,
we get  $\bb\sp{3r + 1 } = \aa\bb\sp{-r} \aa$.
Thus 
$$
\bb\sp r = \aa\bb\sp{-(3r+1)}\aa 
= \aa \bb\sp{-r} \bb\sp{-(2r+1)} \aa = \aa\bb\sp{-r}\aa\sp{-1}.
$$
This implies 
$$
\bb\sp{rq} = \aa\bb\sp{-rq} \aa\sp{-1} 
= \aa(\aa\sp{-2})\sp r \aa\sp{-1} = \aa \sp{-2r} = \bb \sp{-rq}.
$$
Thus $\bb\sp{2rq} = \bb\sp{(q-1)q} = e$.
\qed
\medskip
Consider
the following two matrices in $GL\sb 2 (\C)$;
$$
A =\pmatrix{0 & e(\oover{q}{4rk}) \cr
e(\oover{q}{4rk}) & 0 \cr},
\quad\and\quad
B=\pmatrix{e(\oover{1}{2rk} + \oover{r}{q}) & 0 \cr
0 & e ( \oover{1}{2rk} - \oover{r}{q}) \cr}.
$$
It is easy to check that
$$
A\sp 2 = B\sp q, \quad 
B\sp{qk} = (B\sp{-r} A) \sp{2k}, \quad\and\quad(B\sp{-r} A ) \sp 2 = c I\sb 2
\quad\where\quad c = e(\oover{1}{2rk}).
$$
In particular, $(B\sp{-r} A)\sp 2 $ is contained in the center of 
$GL\sb 2 (\C)$.
Hence there exists a homomorphism 
$$
\vphi : \Hqk \lra GL\sb 2 (\C),
$$
which maps $\aa$ to $A$ and $\bb$ to $B$.
Let
$P : GL\sb 2 (\C) \to PGL\sb 2 (\C)$
be the natural projection.
Then the image of $P\circ \vphi$,
which is generated by
$P(A)$ and $P(B)$,
is the dihedral group $D\sb q$.
Moreover,
$\Ker P \cap \Im \vphi$ contains
a cyclic group of order $2rk$
generated by $(B\sp{-r} A ) \sp 2 = c I\sb 2$.
Hence the order of $\Hqk$ is at least $2rk$ times the order $2q$ of $D\sb q$.
Combining this with Claim,
we see that the order of $\Hqk$ is exactly $2q(q-1)k$,
and $\vphi$ is  injective.
The homomorphism $P\circ \vphi$ gives $\Hqk$ 
a structure of the central extension
of $D\sb q$ by the cyclic group $\Z/(k(q-1))$.
\qed
\bign
{\bf \S 4. Singularities of the curve $C$ in the dihedral case}
\bfclaim{Proposition 3.}
{\sl 
Suppose that $p=m=2$.
Then the singular locus of $C$ consists of
$(2ql-3k)l$ points of type $A\sb{q-1}$.
}
\medn
{\it Proof.}
\hs
Let $\tt$ be a complex parameter,
and we denote by $C\sb\tt\st \Pt$
the curve of degree $2(ql-k)$ defined by
$F(S, T, X, \tt Y)=0$.
Thus the curve $C$ we have been treating  is $C\sb 1$.
We prove Proposition 3 in  the following way;
first, we investigate the singular points of $C\sb 0$, and
then we see what happens to these singular points 
when $\tt$ becomes a non-zero small number $\vee$.
\medskip
The curve $C\sb 0$ is defined by the equation
$$
X (S\sp 2 X + 2 T \sp q) = 0.
$$
It has two  irreducible components.
We put
$$C\o (1) = \{ X = 0 \}, \qquad\and\qquad
C\o(-1) = \{ S\sp 2 X + 2 T\sp q = 0 \}.
$$
Consider the linear pencil
$$
D\sb{\ll} = \{S\sp 2 X + \ll T\sp q = 0 \} \sb {\ll \in \P\sp 1}
$$
of curves of degree $ql$ spanned by the curves 
$\locus{S\sp 2 X = 0}$ and $\locus{T\sp q =0}$.
We have $C\o (-1) = D\sb 2$.
The base locus of this pencil
$$
\Bs = \{ S=T=0 \} \cup \{ X=T=0 \}
$$
consists of $kl + (ql-2k) l$ distinct points.
(Recall that $S$, $T$ and $X$ are general.)
By Bertini's theorem, 
except for finite values $0$, $\infty$, $\ll\sb1$, \dots , $\ll\sb M$ of $\ll$,
the members $D\sb {\ll}$ are non-singular outside $\Bs$.
Since $T$ is general,
we may assume that none of $\ll\sb 1$, \dots, $\ll\sb M$ coincides with 
$2$.
Hence $C\o (-1)$ is non-singular outside $\Bs$.
The intersection points of $C\o(-1)$ and  $C\o(1)$ are also contained in $\Bs$.
Since $X$ is general, $C\o (1)$ is non-singular.
Hence we have $\Sing C\o \st \Bs$.
On the other hand,
it is easy to see that $C\o$ is singular at each point of $\Bs$.
Thus we get
$$
\Sing C\o = \{ S=T=0 \} \cup \{ X=T=0 \}.
$$
\par
Let $P \in \Sing C\o$ be an intersection point of the curves
$\locus{S=0}$ and $\locus{T=0}$.
\bfclaim{Claim 1.}
{\sl
There exist a small open neighborhood $\DD \st \Pt$ of $P$ 
and a small positive real number $r$
such that $C\sbee$ has only one singular point in $\DD$ for any $\vee \in \C$
satisfying $|\vee | < r$.
Moreover, the singular point is of type $A\sb {q-1}$.
}
\medn
{\it Proof.}
\hs
Since $S$, $T$, $X$  and $Y$ are chosen generally, the curves 
$\locus{S=0}$ and $\locus{T=0}$ intersect transversely at $P$,
and both of $X$ and $Y$ are non-zero at $P$.
Let $(u, v)$ be local analytic coordinates of $\Pt$ 
around $P$ such that the curves $\locus{S=0}$ and $\locus{T=0}$
are defined by $u=0$ and $v=0$, respectively.
Then the defining equation of $C\sbee$ 
locally around $P$ is of the form
$$
{{(au\sp 2 + v\sp q)\sp 2 - (\vee b u\sp 2 + v\sp 2)\sp q}\over{u\sp 2}}
= a\sp 2 u\sp 2 + 2a v \sp q - \vee b (u\sp 2 F(u, v, \vee) + q v\sp{2(q-1)} )
 =0,
$$
where $a$ and $b$ are holomorphic functions of $(u,v)$ corresponding to $X$ 
and $Y$,
respectively,
satisfying $a(P)\not = 0$ and $b(P)\not = 0$,
and $F(u, v, \vee)$ is a holomorphic function of
$(u,v, \vee)$ such that $F(0, 0, 0)=0$.
Since $a(P)\not = 0$,
there exist holomorphic functions
$\phi (u, v, \vee)$ and $\psi (u, v, \vee)$,
defined locally around $(u, v, \vee) = (0,0,0)$,
such that
$$
\eqalign{
\phi(u, v, \vee)\sp 2 & = a(u, v) \sp 2 - \vee b (u, v) F(u,v,\vee), \quad 
\and \cr
\psi(u, v, \vee) \sp q & = -2 a(u, v) + \vee q b(u, v) v \sp{q-2}.
}
$$
We put
$$ 
\wt u = u\phi (u,v, \vee) \quad\and\quad \wt v = v\psi(u, v, \vee).
$$
Then $C\sbee$ is defined by
$$
\wt u\hsa\sp 2 - \wt v \hsa\sp q = 0.
$$
Note that
$$
\def\pa#1#2{{\oover{\partial \wt #1}{\partial #2}}(0,0,0)}
\hbox{det} 
\varpmatrix{30pt}{\pa{u}{u} & \pa{u}{v} \cr \pa{v}{u} & \pa{v}{v}}
 \quad\not =\quad 0,
$$
\medn
because of $\phi(0,0,0) \not= 0$ and $\psi (0,0,0) \not= 0$.
This implies that
there exist a small open neighborhood $\DD$ of $P$ on $\Pt$
and a small positive real number $r$ such that  $(\wt u , \wt v)$
are local analytic coordinates in $\DD$ when $|\vee | < r$.
Since the only singular point of $C\sbee = \{ \wt u\sp 2 - \wt v \sp q = 0 \}$
in $\DD$ is the origin $P$ and it is of type $A\sb{q-1}$,
we have completed the proof of Claim 1. \qed
\medskip
Let $Q\in \Sing C\o$ be an intersection point 
of the curves $\locus{X=0}$ and $\locus{T=0}$.
\bfclaim{Claim 2.}
{\sl
There exist a small open neighborhood $\DD \st \Pt$
of $Q$ and a small positive real number $r$
such that,
when $\vee$ satisfies $0< |\vee| < r$,
then $C\sbee$  has exactly $2$ singular points 
in $\DD$, and each of them is of type $A\sb{q-1}$.
}
\medn
{\it Proof.}
\hs
Since $S$, $T$ $X$  and $Y$ are chosen generally,
the curves $\locus{X=0}$ and $\locus{T=0}$ intersect transversely at $Q$,
and both of $S$ and $Y$ are non-zero at $Q$.
Let $(u, v)$ be local analytic coordinates of $\Pt$ 
around $Q$ such that the curves $\locus{X=0}$ and $\locus{T=0}$
are defined by $u=0$ and $v=0$, respectively.
Since $S$ is non-zero at $Q$,
the defining equation of $C\sbee$ 
locally around $Q$ is of the form
$$
(au + v\sp q)\sp 2 - (\vee b + v\sp 2)\sp q = 0,
$$
where $a$ and $b$ are holomorphic functions
corresponding to $S\sp 2$ and $S\sp 2 Y$, respectively, such that
$a(Q) \not = 0$ and $b(Q) \not = 0$.
Since $b(Q) \not = 0$,
there exists a holomorphic function $b\sp{1/2}$ defined
locally around $Q$.
We put
$$
\wt u  = {{au + v\sp q}\over{b\sp{q/2}}}
\quad\and \quad \wt v = {{v}\over{b\sp{1/2}}}.
$$
Then there is an open neighborhood $\DD \st \Pt$
of $Q$ such that $(\wt u, \wt v)$ is a local analytic coordinate system on 
$\DD$.
In $\DD$, 
$C\sbee$ is defined by
$$
\wt u\sp{\hskip 1pt 2} - (\vee + \wt v\sp{\hskip 1pt  2} ) \sp q = 0.
$$
Hence,
if $\vee$ is small enough and non-zero,
the singular points of $C\sbee$ in $\DD$ consists of the two  points
$$
(\wt u , \wt v ) = (0, \sqrt{-\vee}) \quad \and \quad (0, -\sqrt{-\vee}),
$$
and both of them are of type $A\sb{q-1}$.
\qed
\medskip
By Claims 1 and 2,
all singular points  of $C\sbee$ are of type $A\sb{q-1}$
when $\vee$ is non-zero and small enough.
The number of them is 
$$ 
\deg S \cdot \deg T + 2 \deg X \cdot \deg T = (2ql-3k) l.
$$
The locus of all $\tt\in \C$ such that $\Sing C\sb{\tt}$ consists of 
$(2ql-3k) l$ points of type $A\sb{q-1}$
is then  Zariski open dense in $\C$.
Since $S$, $T$ $X$  and $Y$ are chosen generally,
we  can conclude that $\tt =1$ is contained in this locus.
\qed
\medskip
Now we can compute the geometric genus $g$ of the curve $C$
in Proposition 3. It is given by
$$
g= 1 +3k+2k\sp 2 - 3l -4kl+2l\sp 2 -6lr - 5klr+6l\sp 2 r + 4l\sp2 r\sp 2,
$$
where $r=(q-1)/2$.
In particular, the geometric genus of the curve $C(q, k)$
in Theorem is
given by
$$
g= 1 - 6kr + k\sp 2 r + 4 k\sp 2 r\sp 2.
$$
Hence, except for the case $(q, k)=(3,1)$ of the three cuspidal plane quartic,
the curve $C(q, k)$ is not rational.

\bign
\centerline{\bf References}
\medskip
\item{\Abhyankar}
{\it S.\ Abhyankar},
Tame coverings and fundamental groups of algebraic varieties,
Part V : Three cuspidal plane quartic,
Amer.\ J.\ Math.\ {\bf 82} (1960), 365 - 373.
\item{\BruceGiblin}
{\it J.\ W.\ Bruce, P.\ J.\ Giblin},
   A stratification of the space of plane quartic curves,
Proc.\ London Math.\ Soc. {\bf 42} (1981), 270-298.
\item{\Degtyarev}
{\it A.\ Degtyarev},
Quintics in $\C \P \sp 2$ with nonabelian fundamental group,
preprint.
\item{\Dimca}
{\it A.\ Dimca},
Singularities and Topology of Hypersurfaces,
Springer-Verlag, Berlin, 1992.
\item{\Nemethi}
{\it A.\ N\'emethi}, 
   On the fundamental group of the complement
of certain singular plane curves,
Math.\ Proc.\ Cambridge Philos.\ Soc.\
{\bf 102} (1987), 453 - 457.
\item{\Nori}
{\it M.\ Nori},
Zariski's conjecture and related problems,
Ann.\ Sci.\ \'Ecole.\ Norm.\ Sup. {\bf 16} (1983), 305 - 344.
\item{\Oka}
{\it M.\ Oka},
   Some plane curves whose complements
have non-abelian fundamental groups,
Math.\ Ann.\  {\bf 218} (1975), 55 - 65.
\item{\Okapre}
{\it M.\ Oka},
Two transformations of plane curves and their fundamental
groups,  
preprint.
\item{\ShimadaF}
{\it I.\ Shimada},  
   Fundamental groups of open algebraic varieties,
to appear in Topology.
\item{\ShimadaW}
{\it I. Shimada},
   A weighted version of Zariski's hyperplane section theorem
and the fundamental groups of complements of plane curves,
MPI preprint 95-26, eprint 9505004.
\item{\Zariski}
{\it O.\ Zariski},
   On the problem of existence of 
algebraic functions of two variables
possessing a given branch curve, 
Amer.\ J.\ Math.\  {\bf 51} (1929), 305 - 328.
\bign
Max-Planck-Institut f\"ur Mathematik
\parn
Gittfried-Claren-Str. 26
\parn 
53225 Bonn, Germany
\parn
shimada@mpim-bonn.mpg.de
\bign
{\it   home address}
\parn
Department of Mathematics
\parn
Faculty of Science
\parn
Hokkaido University
\parn
Sapporo 060, Japan
\parn
shimada@math.hokudai.ac.jp
\end